\documentclass[11pt]{article}         
\usepackage{amsfonts,amsmath,amssymb}  

\title{Symmetry for solutions of two-phase semilinear elliptic equations on hyperbolic space}

\author{Isabeau Birindelli \thanks{Partially supported by PRIN project - Viscosity, Metric 
and Control Theoretic Methods in Nonlinear PDE (MIUR, Italy)} \\ Sapienza Universit\`a di Roma 
\and Rafe Mazzeo \thanks{Supported by the NSF through grant DMS-0505709} 
\\ Stanford University}

\date{}

\newtheorem{theo}{Theorem}[section]
\newtheorem{prop}[theo]{Proposition}

\newtheorem{cor}[theo]{Corollary}
\newtheorem{lemme}[theo]{Lemma}

\newcommand{\R}{{\rm I}\!{\rm  R}}
\newcommand{\del}{\partial}
\newcommand{\calC}{\mathcal C}
\newcommand{\HH}{\mathbb H}
\newcommand{\RR}{\mathbb R}
\newcommand{\e}{\epsilon}
\newcommand{\TT}{\mathbb T}
\def\R{{\rm I}\!{\rm  R}}
\def\la1{\lambda_1}

\def\tw{\tilde w}
\def\tI{\tilde{I}}
\def\tu{\tilde{u}}
\def\calU{\mathcal U}
\def\calV{\mathcal V}

\begin{document}
\maketitle

\begin{abstract}
Assume that $f(s) = F'(s)$ where $F$ is a double-well potential. Under certain
conditions on the Lipschitz constant of $f$ on $[-1,1]$, we prove that 
arbitrary bounded global solutions of the semilinear equation $\Delta u = f(u)$ 
on hyperbolic space $\HH^n$ must reduce to functions of
one variable provided they admit asymptotic boundary values on $S^{n-1} = 
\del_{\infty}\HH^n$ which are invariant under a cohomogeneity one subgroup
of the group of isometries of $\HH^n$. We also prove existence of these 
one-dimensional solutions. 
\end{abstract}

\section{Introduction}
The equation $-\Delta u + u(u^2-1) = 0$ is the Euler-Lagrange equation for the 
functional $u \mapsto \int \frac12 |\nabla u|^2 + \frac14 (u^2-1)^2$; it
has been intensively studied for many years, and has many applications in geometry, 
physics and elsewhere. The important feature of this energy functional for many 
of these investigations is the `double-well' behaviour of the nonlinear term 
$F(u) = (1/4)(u^2-1)^2$, and indeed, many results generalize directly to extremals 
of other energy functionals where $F$ is replaced by a more general function with 
qualitatively similar behaviour.

One important theme in this research is to study the symmetry properties of solutions, 
in particular those which are defined globally on $\RR^n$. Most famous here is the De Giorgi 
conjecture, which asserts that, so long as $n \leq 8$, any solution $u(x)$ of this PDE 
which takes values in $[-1,1]$, and satisfies certain monotonicity properties to be made explicit 
below, necessarily has the form $U(\ell \cdot z)$, where $\ell \in \RR^n$ is constant, $z \in \RR^n$ 
and $U(s)$ satisfies $U'' = F'(U)$; in  other words, $u$ reduces to a function of one variable. 
This has been proved for $n \leq 3$ \cite{gg,ac}. 
The cases $4\leq n<8$ have been proved by Savin \cite{S}, but with an additional assumption 
beyond monotonicity that the solution converges to $\pm 1$ at infinity. For related results 
see \cite{aac,f,bhm}, and \cite{fv} contains a survey. 

Our goal here is to investigate these symmetry results from a broader perspective.
Let $(X,g)$ be a complete homogeneous Riemannian manifold, and $G = \mbox{Isom}\,(X)$ 
its group of isometries (which thus acts transitively). Consider bounded solutions $u$
of the equation
\begin{equation}
\Delta_g u - F'(u) = 0,
\label{eq:heq}
\end{equation}
where $F$ is a double-well potential function. More precisely, 
\begin{equation}
\begin{array}{rl}
& F \in \calC^2(\RR),\ F \geq 0,\ F(s) \to \infty\ \mbox{when}\ |s| \to \infty, \\
& F^{-1}(0) = \{-1,+1\}\ \mbox{and}\ F''(\pm 1) > 0.
\end{array}
\label{eq:condF}
\end{equation} 
We shall also always assume that $s F'(s) \geq 0$ when $|s| \geq 1$, which is
a technical convenience which (by the maximum principle) allows us to assume only 
that $u$ is bounded rather than that $u(z) \in (-1,1)$ for all $z$; furthermore, 
we also assume that the only other critical point of $F$ is a local maximum at 
$s=0$, though this too can be weakened. The main question is as follows:
\begin{quote}
If $u$ is a bounded solution of (\ref{eq:heq}) defined on all of $X$, then under what
conditions is $u$ necessarily invariant under some subgroup of positive dimension of $G$?
In particular, under what conditions is $u$ invariant under a subgroup $H \subset G$ which 
has orbits of codimension one in $X$?
\end{quote}
The constant solutions $u \equiv \pm 1$ are invariant under $G$ itself, but we are more
interested in solutions which `jump states', i.e.\ such that $u \to -1$ in certain directions 
at infinity and $u \to +1$ in others.

The types of hypotheses we shall place on $u$ concern its `asymptotic boundary values', i.e.\ 
its limits in different directions at infinity. This differs from the hypotheses in the De Giorgi 
conjecture where one assumes only that $u$ is monotone along some fixed family of parallel lines;
our hypothesis is analogous to the so called Gibbons condition \cite{bhm}. Part of the work
below involves showing that once one has imposed certain asymptotic boundary values, the solution
necessarily has a monotonicity property; however, as we indicate below, for certain choices of 
boundary values, if $F''(0)$ is sufficiently large, there are ODE solutions of (\ref{eq:heq}) 
which are not monotone but do have these correct asymptotic boundary values.

We shall not formulate the precise problems we study in this generality, but consider only
the particular case that $X$ is the $n$-dimensional hyperbolic space $\HH^n$, and 
hence $G = \mbox{Isom}\,(\HH^n) = \mbox{O}(n,1)$. This is of interest since the isometry group 
$G$ is substantially more complicated than the group $E(n)$ of Euclidean motions. Many of the
results here can be proved in much greater generality; for example, Theorem 2
below remains valid for any rank one symmetric space $X$ (and could even be formulated suitably
for higher rank symmetric spaces of noncompact type). The key property needed is that $G$ contains a sufficiently
large collection of reflections. However, for simplicity we focus on this problem on hyperbolic space. 

The geodesic compactification $\overline{\HH^n} = \HH^n \cup S^{n-1}$ is obtained by adding
a point at the end of each geodesic ray emanating from some fixed point. This asymptotic
boundary has a natural topology and conformal structure, and isometries of 
$\HH^n$ induce conformal transformations of the sphere at infinity and conversely. 
A function $u$ on $\HH^n$ has asymptotic boundary value $\phi$ on $S^{n-1}$ if $u$ 
is continuous (a.e.) on $\overline{\HH^n}$ and its restriction to the boundary is equal to $\phi$. 

Briefly, we consider solutions $u$ of (\ref{eq:heq}) on $\HH^n$ with asymptotic boundary values 
invariant with respect to a cohomogeneity one subgroup $H \subset G$, and the main question we ask 
is whether the solution $u$ itself is also invariant with respect to $H$. Since $\HH^{n}/H$ is 
one-dimensional, this means that $u$ is determined by a function of one variable (which hence 
satisfies an ODE). 

Before stating our results precisely, we recall the three main models for $\HH^n$: the upper-half-space,
\[
\RR^n_+ = \{z = (x,y) \in \RR^n: x > 0,\ y \in \RR^{n-1}\}, \qquad g = \frac{dx^2 + |dy|^2}{x^2},
\]
the Poincar\'e ball, 
\[
B^n = \{z \in \RR^n: |z| < 1\}, \qquad g = \frac{4|dz|^2}{(1-|z|^2)^2},
\]
and via geodesic normal coordinates, 
\[
\RR^n = \{z = r\omega: r \geq 0,\ \omega \in S^{n-1}\}, \qquad g = dr^2 + \sinh^2 r\, d\omega^2.
\]
There are three subgroups of $G$ which have orbits in $\HH^n$ of codimension one:
\[
H_e = \mbox{O(n)}, \qquad H_p = E(n), \qquad H_h = \mbox{O}(n-1,1).
\]
The subscripts `e', `p' and `h' stand for elliptic, parabolic and hyperbolic, respectively. 
The action of $H_e$ is by rotations around a fixed point $p \in \HH^n$ or by reflections across a 
totally geodesic hyperplane containing that point; in the ball model, with $p$ at the origin, or in
geodesic polar coordinates, this corresponds to the standard rotation group acting on $B^n$ or $\RR^n$,
respectively. Next, $H_p$ acts by Euclidean motions on each horosphere centered
at a point $q \in S^{n-1} = \del \HH^n$; in the upper-half-space model, with $q$ the point at infinity, these
horospheres are the hyperplanes $x = \mbox{constant}$ and the action appears as the standard linear one in
each $\RR^{n-1}$. Finally, the inclusion of $H_h$ in $G$ depends on the choice of a totally geodesic hyperplane
$P \cong \HH^{n-1} \subset \HH^n$. Any isometry $h$ of $P$ extends (uniquely up to a reflection across $P$) to
an isometry of the entire space $\HH^n$, and the subgroup $H_h$ consists of all isometries obtained in this way.
In the ball model, we could take $P = \{z_n = 0\}$, while in the upper half-space model a convenient choice is 
$P = \{|(x,y)| = 1\}$. The hypersphere $P \cap \del \HH^n = S^{n-2} \subset S^{n-1}$ divides the boundary sphere 
into two components, which we denote by $S_-$ and $S_+$, respectively. Finally, we let $t$ be the function
on $\HH^n$ which gives the signed distance from $P$, so $t > 0$ on one component of $\HH^n \setminus P$ and
$t < 0$ on the other.

Any function on $S^{n-1} = \del \HH^n$ invariant under the induced action of $H_e$ must be constant; a function 
invariant with respect to the action of $H_p$ can assume only two values, one at the parabolic fixed point
$q$ and the other on the other orbit, $S^{n-1}\setminus \{q\}$; finally, any function on $S^{n-1}$ which is
invariant under the subgroup of $H_h$ which does not interchange the two components of $\HH^n \setminus P$ 
takes on one constant value on all of $S_-$ and another on all of $S_+$. 
We present two main theorems. The first guarantees the existence of $H_i$ invariant solutions for each of 
these subgroups, $i = e, p, h$, and the second gives conditions under which all solutions (with $H_i$-invariant 
boundary values) must themselves be $H_i$ invariant. This second theorem includes some decay hypotheses
on the solution, and to show that these are necessary we also prove the existence of non-symmetric
solutions which have just slightly weaker decay. 

Before stating these, we introduce a few important constants associated to this problem. Let $L$
denote the Lipschitz constant of $f$ on the interval $[-1,1]$, and $\lambda = -f'(0) > 0$. Consider
the two equations:
\begin{equation}
\alpha^2 - (n-1)\alpha + \lambda = 0, \qquad \alpha^2 - (n-1)\alpha + L = 0,
\label{eq:quadratics}
\end{equation}
and denote the roots of these by $\alpha_\pm$ and $\beta_\pm$, respectively. Note that it is
always true that $\lambda \leq L$. We shall need to assume at various points below that
$L \leq (n-1)^2/4$; this ensures that all of these roots are real, and satisfy
\begin{equation}
0 < \alpha_- \leq \beta_- \leq \frac{n-1}{2} \leq \beta_+ \leq \alpha_+ < n-1.
\label{eq:indroots}
\end{equation}

\begin{theo} Consider solutions of (\ref{eq:heq}) which are invariant under $H_e$, $H_p$ or $H_h$,
where $f = F'$ and $F$ satisfies the hypotheses in (\ref{eq:condF}). 
\begin{itemize}
\item[i)] When $i=e$, the only ODE (i.e.\ radial) solution which has constant asymptotic boundary value $\pm 1$
is the constant solution $u \equiv \pm 1$. There is a one-parameter family of radial solutions, all
of which decay like $(1-|z|)^{\alpha_-}$ as $|z| \to 1$ in the ball model; however, if $|u| \leq C(1-|z|)^\alpha$ 
for some $\alpha > \alpha_-$, then $u \equiv 0$.
\item[ii)] When $i=p$, there exists an $H_p$-invariant solution $u_p$ which is unique up to translation,
and which tends to $0$ on $S^{n-1}\setminus \{q\}$ and to $+1$ at $\{q\}$, where $q$ is the parabolic fixed
point of $H_p$. For generic functions $f$, this solution decays like $x^{\alpha_-}$ as $x \searrow 0$. 
On the other hand, there exists no $H_p$-invariant solution which tends to $-1$ along $\RR^{n-1}$ 
and $+1$ at infinity, or vice versa.
\item[iii)] When $i=h$, there exists an $H_h$-invariant solution $u_h$ which is monotone as a function of
the signed distance $t$, and which tends to $-1$ in $S_-$ and $+1$ in $S_+$. 
\end{itemize}
\end{theo}
Note that we are not asserting the uniqueness of the ODE solution in case iii), and indeed it is likely
that there is a sequence of `multibump' solutions of this problem which have an arbitrary
finite number of oscillations. 

The converse statement is provided by the
\begin{theo} 
Let $u$ be an arbitrary bounded solution of (\ref{eq:heq}), where $f = F'$ and $F$ satisfies the hypotheses in 
(\ref{eq:condF}). 
\begin{itemize}
\item[i)] Let $\eta = -1$ or $+1$. If $|u(z) - \eta| \to 0$ as $z \to \infty$ in $\HH^n$, then 
$u \equiv \eta$. On the other hand, suppose that $f$ is Lipschitz with Lipschitz constant $L 
\leq (n-1)^2/4$ and $f(s) \in \calC^{1,\delta}$ in some neighbourhood of $s=0$ for some 
$\delta > 0$. If $|u(z)|/(1-|z|)^\alpha \leq C$ for some $\alpha > \alpha_-$, where $z$ is the 
coordinate in the ball model), then $u \equiv 0$. It suffices to require only that $f$ is Lipschitz
if this exponent of decay satisfies $\alpha > \beta_-$. 
\item[ii)] Suppose that, as before, $f(t)$ is Lipschitz with Lipschitz constant $L \leq (n-1)^2/4$.
If $|u| \leq C x^\alpha$ on all of $\HH^n$ for some $\alpha > \alpha_-$ (where we use the 
upper half-space variable $x$), and $u \to +1$ in any conic region $x > C_1|y - y_0| + 
C_2$ as $(x,y) \to \infty$, then $u$ is invariant under $H_p$, and hence is a function of $x$ alone. 
\item[iii)] If $u \to -1$ uniformly on compact subsets of $S_-$ and $u \to +1$ on
compact subsets of $S_+$, then $u$ depends only on the signed distance $t$ from
the totally geodesic $\HH^{n-1}$. 
\end{itemize}
\end{theo}

If $u$ is an arbitrary solution which is bounded by $x^\alpha$ for $\alpha > \alpha_-$, then by
part ii), $u$ is invariant under $H_p$; however, by part ii) of the preceding theorem, this
ODE solution does not decay so quickly as $x \searrow 0$, at least for most choices of 
function $f$. The conclusion is simply that there are no solutions with this rate of decay.

As was noted earlier, unlike for the De Giorgi conjecture in $\RR^n$, we assume conditions about 
the asymptotic boundary values but do not impose monotonicity conditions on the solution. 

Our final result is to show that the decay hypotheses in parts i) and ii) of Theorem 1.2 are necessary.
\begin{theo}
Let $f$ be as above.
\begin{itemize}
\item[i)] There is an infinite dimensional family of solutions of (\ref{eq:heq}), each element of which
decays like $(1-|z|)^{\alpha_-}$ as $|z| \to 1$ in the ball model of $\HH^n$.
\item[ii)] There is an infinite dimensional family of solutions of (\ref{eq:heq}), each element of which
is bounded by $x^{\alpha_-}$ in the upper half-space model of $\HH^n$.
\end{itemize}
\end{theo}

The equation (\ref{eq:heq}) has many other interesting global bounded solutions on $\HH^n$. 
It would be interesting to construct bounded solutions which are piecewise constant on the 
sphere infinity, taking only the values $0, \pm 1$ there. The results here address only 
the very simplest situations of this sort. It would also be interesting to understand 
if monotonicity but no assumptions on the asymptotic boundary conditions are enough to
imply that the solution is one-dimensional.  We shall come back to these issues elsewhere.

\section{ODE solutions}
We now prove the existence of one-dimensional solutions to (\ref{eq:heq}) which
are invariant under either $H_p$ or $H_h$, and discuss the uniqueness of solutions 
in the $H_e$ case.

\subsection{The elliptic case}
The existence of radial $H_e$ invariant solutions is obvious, since we can simply let
$u \equiv 0$ or $\pm 1$. We now consider the uniqueness of these ODE solutions.

\begin{prop}
The only $H_e$ solution of (\ref{eq:heq}) which tends to $\pm 1$ at infinity is the
constant solution. However, there is a one-parameter family of radial solutions to
this equation which tend to $0$ at infinity.
\end{prop}

\noindent{\bf Proof:}
For radial functions and in geodesic normal coordinates, (\ref{eq:heq}) takes the form
\begin{equation}
u'' + (n-1)\coth r \, u' = f(u).
\label{eq:eqradial}
\end{equation}
Suppose that $u(r) \to \pm 1$ as $r \to \infty$. Multiplying by $u'$ and integrating from $0$ to $A$ yields
%\begin{equation*}
\begin{multline*}
\frac12 u'(A)^2 + (n-1)\int_0^A \coth r\, u'(r)^2\, dr \\ = \int_0^A f(u(r))u'(r)\, dr = F(u(A)) - F(u(0)),
\end{multline*}
%\end{equation*}
where we used that $u'(0) = 0$. As $A \to \infty$, $u'(A) \to 0$ and $F(u(A)) \to 0$, whereas
$F(u(0)) \geq 0$, so taking the limit gives
\[
(n-1) \int_0^\infty \coth r\, u'(r)^2\, dr = -F(u(0)) \leq 0;
\]
this implies that $u' \equiv 0$, hence $u \equiv \pm 1$. 

This calculation does not provide any information when $u \to 0$ at infinity since $\lim F(u(A)) > 0$.
In fact, there is a one-parameter family of nonconstant radial solutions tending to $0$ as $r \to \infty$. 
To see this, consider the solution $u_a(r)$ of (\ref{eq:eqradial}) which satisfies the initial
conditions $u_a(0) = a \in (-1,1)$, $u_a'(0) = 0$. By the maximum principle, $|u_a(r)| < 1$ for
all values of $r$ for which the solution exists, and because of this boundedness it is easy to
see that $u_a'$ and $u_a''$ also remain bounded, so that $u_a(r)$ exists for all $r \in \RR^+$. 

We next show that $u_a$ has a limit as $r \to \infty$. If it were monotone there would be nothing to prove otherwise there would exist a sequence $\{r_j\}$ tending to infinity where $u_a'$ vanishes. There
would in fact exist some subsequence $r_j'$ such that $u_a(r_j')$ tends to a limit (and $u_a'(r_j') = 0$
for all $j$). Multiply (\ref{eq:eqradial}) by $u_a'$ and perform exactly the same integration by parts as 
above on the interval $(r_j',r_{j+1}')$ to get
\[
(n-1)\int_{r_j'}^{r_{j+1}'} \coth r\, u'_a(r)^2\, dr = F(u_a(r_{j+1}')) - F(u_a(r_j')) \to 0.
\]
Since the left hand side is nonnegative, the sequence $\{F(u_a(r_j'))\}$ is monotone increasing.
If the lengths of all the intervals $|r_{j+1}' -r_j'| \to 0$, then since $u_a'$ is bounded, we
conclude that $u_a$ itself has a limit. Otherwise, we conclude at least that the $L^2$ norm
of $u_a'$ is uniformly small on each such interval, and hence it converges uniformly to $0$  too, here too we conclude that $u_a$ converges.

Finally to prove that $u_a$ converges to zero, writing (\ref{eq:eqradial}) 
\[ [g(r)u'_a]'=g(r)f(u_a)\]
with $g(r)=(\sinh r)^{n-1}$
and integrating it on each $(r_j',r_{j+1}')$ gives
\[0=\int_{r_j'}^{r_{j+1}'}g(r)f(u_a(r))dr.\]
This implies that for each $j$ there exists $r_{o,j}\in (r_j',r_{j+1}')$  such that $f(u_a(r_{o,j}))=0$. This gives that $u_a$ converges to either $0$ or $\pm 1$, and since we already proved that it cannot converge to $\pm 1$
without being identically equal to $\pm 1$, we conclude that $u_a(r) \to 0$ as $r \to \infty$,
as claimed.  \hfill $\Box$ 

\subsection{The parabolic case}
Consider $\HH^n$ in the upper half-space model, placing the parabolic fixed point at infinity. 
An $H_p$-invariant solution satisfies the equation
\begin{equation} 
\label{ham0} 
x^2 u''(x) + (2-n)x u'(x) - f(u) = 0,
\end{equation}
where $x \in \RR^+$. We assume that $f$ satisfies the conditions in (\ref{eq:condF}), and moreover 
that $f'(0) < 0$. 

\begin{prop}
There are no bounded solutions of (\ref{ham0}) which satisfy 
\[
\lim_{x\to 0}u(x)=-1,\  \lim_{x\to +\infty}u(x)=1. 
\]
On the other hand there exists a solution, unique up to the rescaling $x \mapsto cx$, $c > 0$, such that
\[
\lim_{x\to 0}u(x)=0,\  \lim_{x\to +\infty}u(x)=1. 
\]
\label{pr:pc}
\end{prop}

\noindent{\bf Proof:} Setting $\xi =\log x$, and letting primes now refer to
derivatives with respect to $\xi$, then (\ref{ham0}) becomes
\begin{equation}
u''(\xi) - (n-1)u'(\xi) - f(u) = 0.
\label{ham01}
\end{equation}
Setting $u' = v$ yields the equivalent autonomous system
\begin{equation} 
\label{ham} 
\begin{array}{l} u'=v\\
v'=(n-1)v+f(u),
\end{array}
\end{equation}
which we shall analyze by phase plane methods. There are three critical points: $(0,0)$,
which is either an unstable node or an unstable spiral depending on whether $-f'(0)$ is less 
than or equal to, or greater than, $(n-1)^2/4$, and $(\pm 1,0)$, which are both saddle points 
with unstable eigenvectors of the form $(1,\alpha_\pm)$, where both $\alpha_\pm > 0$. 
We are thus looking for heteroclinic orbits connecting any two of these critical points.

Let us first examine the integral curve which emanates from $(-1,0)$ as $\xi \to -\infty$ and
moves into the upper half-plane $v > 0$. There are only three possible scenarios: either this 
orbit crosses the line $v=0$ in the open interval where $0 < u < 1$, or it tends to $(1,0)$ 
as $\xi \to +\infty$, or finally that it passes above this point. Note that the orbit cannot cross
in the interval $-1 < u < 0$ since the vector field points upward there. We shall show that the
first two possibilities cannot occur, which then proves the first assertion of the theorem,
as well as the second since it then follows that precisely one of the integral curves emanating 
from $(0,0)$ as $\xi \to -\infty$ converges along the stable direction to $(1,0)$, and this is 
the solution we want. Note that as $-f'(0)$ increases, this solution may spiral around the
origin more and more times; indeed, once $-f'(0) > (n-1)^2/4$, it spirals around infinitely
often. There is a critical value $\gamma$ so that if $-f'(0) \leq \gamma$, then
the corresponding solution stays strictly positive, but if $-f'(0) > \gamma$, then the solution 
becomes negative;  hence one should not expect solutions to be monotone unless $|f'(0)|$ is sufficiently small. 

Thus suppose that $(u(\xi),v(\xi))$ is a solution curve connecting $(-1,0)$ to $(+1,0)$. 
Multiply (\ref{ham01}) by $u'$ and integrate from $-A$ to $+A$ to get
\[
\begin{split}
\frac12 (u'(A)^2 - u'(-A)^2) - (n-1)\int_{-A}^A (u')^2\, d\xi = \\
\int_{-A}^A f(u(\xi))u'(\xi)\, d\xi = F(u(A)) - F(u(-A)).
\end{split}
\]
Taking the limit as $A \to \infty$ yields that
\[
\int_{-\infty}^\infty (u')^2 = 0,
\]
since $F(-1) = F(+1)$. This proves that such a solution cannot exist.  If the integral curve 
emanating from $(-1,0)$ reaches the $u$-axis in the range $0<u<1$ at some finite time $T$, 
then since $u'(T) = v(T) = 0$, the same calculation yields that
\[
-(n-1)\int_{-\infty}^T (u')^2\, d\xi  = F(u(T)) - F(-1) > 0,
\]
which is again impossible. The same reasoning shows that there cannot be a homoclinic orbit connecting 
$(1,0)$ to itself. All of this establishes the existence of a unique integral curve
connecting $(0,0)$ to $(1,0)$.

Unless $f$ is very special, this solution will not approach $(0,0)$ along the eigenvector
with eigenvalue $\alpha_+$, hence for generic $f$, $u(x) \sim c x^{\alpha_-}$ as $x \searrow 0$. 

The case of a solution going from $(0,0)$ to $(-1,0)$ is handled similarly.  \hfill $\Box$

\medskip

There is one specific case where this solution can be written explicitly. Let $f(u)=ku(1-u^2)$, 
where $k=\frac{2(n-1)^2}{9}$, and set $a=\frac{n-1}{3}$; then 
\[
u(x)=\frac{x^a}{1+x^{a}}
\]
is an explicit solution of (\ref{ham0}) satisfying
\[
\lim_{x\rightarrow 0}u=0,\quad \lim_{x\rightarrow \infty}u=1.
\]
(Starting the `time' variable $\xi$ at different points along this solution yields the 
family of solutions $x^a/(b + x^a)$ for any $b > 0$.) 

\begin{cor} If $u$ is any bounded solution of (\ref{eq:heq}) such that, in the upper half-space model,
\[
\lim_{x\to 0}u=-1, \quad  \lim_{x\to \infty}u=1,
\]
then $u$ does not reduce to a function of only one variable.
\end{cor}

\subsection{The hyperbolic case} The solutions of (\ref{eq:heq}) which are invariant under the 
subgroup $H_h$ satisfy the ODE 
\begin{equation}\label{hyp}
U''(t) + (n-1)\tanh t \, U' = f(U), \qquad \lim_{t \to \pm \infty} U(t) = \pm 1.
\end{equation}
We now prove the existence of such a solution. 

\begin{prop}
Suppose that $f(s) = F'(s)$, where $F$ satisfies the same assumptions as before. Then there exists a 
solution $U$ of (\ref{hyp}).
\label{pr:hc}
\end{prop} 
\noindent {\bf {Proof.}} We first find a function $U_T(t)$ defined on the interval $[-T,T]$ for
any  $T > 0$ which minimizes the functional
\[
E(u) = \int_{-T}^T( \frac12 |u'(t)|^2 + F(u(t))\, )(\cosh t)^{n-1}\, dt
\]
over the class of functions $u \in H^1([-T,T]; (\cosh t)^{n-1}\, dt)$ satisfying $u(\pm T) = \pm 1$. 
Since $F$ is nonnegative, it is standard that such a function exists. Note that $U_T(t)$ satisfies the 
ODE in (\ref{hyp}) and $U_T(t) \in (-1,1)$ for $|t| < T$. It is also not hard to see that $U_T$ is 
monotone increasing. Indeed, by the maximum principle, no solution of this ODE can have a positive local 
minimum or negative local maximum. Therefore, if $U_T$ is not monotone, then at least it must
be monotone increasing on two intervals $[-T,A]$ and $[B,T]$, where $U_T \leq 0$ on $[-T,A]$,
$U_T \geq 0$ on $[B,T]$ and $U_T(A) = U_T(B) = 0$. If monotonicity were to fail, it would
be because $U_T$ oscillates on $[A,B]$. However, this cannot be the case for a minimizer,
since the function which agrees with $U_T$ on $[-T,A] \cup [B,T]$ and is identically $0$
on $[A,B]$ has less energy. (There are quite likely nonminimal critical points of this functional
which do oscillate.) 

In order to prove the existence of $\lim_{T \to \infty} U_T$, it suffices to show
that the value $U_T(0)$ remains bounded away from $-1$ and $+1$; alternately, we must show that
the value $a = a_T$ where $U_T(a_T) = 0$ remains bounded as $T \to \infty$. If this were
not the case, consider the sequence of functions $U_T(\cdot - a_T)$; this has a limit
as $T \to \infty$ which is defined either on all of $\RR$ or else on some half-line
$[-A,\infty)$ or $(-\infty,B]$, and has the boundary values or limits $-1$ and $+1$ 
at the left and right ends in each case. To be concrete, suppose that $a_T \to -\infty$,
so that the limit function $v(t)$ satisfies $v'' - (n-1)v' - f(v) = 0$ on 
either $\RR$ or $[-A,\infty)$. Multiply this equation by $v'$ and integrate over the domain 
of definition to get
\[
-\frac12 |v'(-A)|^2 - (n-1)\int_{-A}^\infty (v')^2\, dt = F(1) - F(-1) = 0
\]
(the first term on the left is absent if $A = \infty$), which is a contradiction. 
This proves that the limit $U$ of $U_T$ exists, satisfies the ODE (\ref{hyp}) and lies
in $(-1,1)$ for all $t$. 

It remains to prove that it has the correct asymptotic limits. Since $U(t)$ is monotone
increasing, there is some limiting value $\omega = \lim_{t \to \infty}U(t)$. However, the
eventuality $\omega < 1$ is ruled out by the ODE again since the limit of the left side
would vanish while $\lim_{t \to \infty} f(U(t)) = f(\omega) \neq 0$.   The corresponding statement is true
for the limit of $U$ as $t \to -\infty$. \hfill $\Box$

\medskip

We remark that the key property of $f$ required in this proof is that $\int_{-1}^1 f(s)\, ds = 0$.

\section{One-dimensional symmetry of solutions}
We now turn to the problem of determining when global bounded solutions of (\ref{eq:heq}) are
invariant under one of the subgroups $H_e$, $H_p$ or $H_h$. The first two cases will require an
additional decay hypothesis which we shall demonstrate later is necessary. This can be stated in 
several equivalent ways. The most convenient is in terms of two different families of 
generalized eigenfunctions of the Laplacian: the plane wave solutions $\sigma_\alpha$ and the 
spherical functions $\phi_\alpha$, which we now define. In upper half-space coordinates, the plane wave 
solution with pole at infinity and exponent $\alpha \in {\mathbb C}$ is defined by
\begin{equation}
\sigma_\alpha(x,y) = x^\alpha;
\label{eq:sigma}
\end{equation}
it satisfies $\Delta \sigma_\alpha = \alpha(\alpha - (n-1))\sigma_\alpha$. Thus, when 
$\mbox{Re}\, \alpha > 0$, $\sigma_\alpha \to 0$ when $x \searrow 0$ and $\sigma \to \infty$
when $x \nearrow \infty$. Similarly, for any 
$\alpha \in {\mathbb C}$ there exists a unique radial ($H_e$-invariant) function $\phi_\alpha$
on $\HH^n$ which satisfies this same eigenfunction equation; when $\alpha < \frac12 (n-1)$, 
$\phi_\alpha \sim (1-|z|)^{\alpha}$ as $|z| \to 1$. 

In the ball model, the distance to the boundary sphere in the Euclidean metric is comparable to the 
exponential of the negative of the distance function for the hyperbolic metric:
\[
C_1 (1-|z|) \leq e^{-d(z,0)} \leq C_2(1-|z|),
\]
so it is not hard to express the decay of $\sigma_\alpha$ and $\phi_\alpha$ in terms
of this geodesic distance. 

We first present two simple lemmas which will be used in several places below.
\begin{lemme}
Let $u$ be a bounded solution of (\ref{eq:heq}) on $\HH^n$. Suppose that $T_j$ is any sequence
of M\"obius transformations. Then the functions $u \circ T_j$ are also solutions of (\ref{eq:heq})
and some subsequence of these converges (at least in $\calC^2$) on every compact set to another
solution of this same equation.
\label{le:lemma1}
\end{lemme}
This follows from the boundedness of the sequence $u \circ T_j$ and standard elliptic estimates.

\begin{lemme}
Let $I(z)$ be a bounded continuous function on $\HH^n$ which satisfies $I \geq c > 0$ and 
has the property that for any sequence 
of M\"obius transformations $T_j$, some subsequence of $I \circ T_j$ converges uniformly on compact
sets. Let $X$ be a vector field on $\HH^n$ with the analogous property, that for any such sequence $T_j$,
the sequence of vector fields $X_j$ obtained by pushing forward $X$ by the action of $T_j$ also has
a convergent subsequence. Suppose that $u$ is a solution to $\Delta u + X \cdot \nabla u = 
I(z)u$ which is bounded below. Then $\inf u \geq 0$ and this infimum equals $0$ if 
and only if $u \equiv 0$. The same conclusion holds if $u$ is only a supersolution to this equation, 
but then we must also assume that $u \circ T_j$ has a convergent subsequence for every sequence $T_j$.
\label{le:lemma2}
\end{lemme} 
\noindent{\bf Proof:} Let $z_j$ be a sequence in $\HH^n$ such that $u(z_j) \to \inf u$. Choose a sequence
of M\"obius transformations $T_j$ so that $T_j(0) = z_j$, and consider the sequence $u_j = u \circ T_j$,
which solves the equation $\Delta u_j + X_j \cdot \nabla u_j = I_j u_j$, where $I_j$ and $X_j$ are
the pullbacks of $I$ and $X$ by $T_j$. Using the hypothesis, choose
a subsequence so that $I_j$ converges locally uniformly to a continuous nonnegative function $\bar{I}$. 
In some fixed neighbourhood of $0$, $u_j$ is also bounded above, and by elliptic estimates we may assume that 
$u_j$ also converges in this neighbourhood to some function $\bar{u}$. Now $\Delta \bar{u} = \bar{I}\, \bar{u}$,
$\bar{I} \geq c>0$ and $\bar{u}$ attains a negative minimum at $0$, which is impossible. 
\hfill $\Box$ 

\subsection{The elliptic case}
The first symmetry result is proved by a `stationary plane' reflection argument, not too different from
the standard moving plane method in the Euclidean case \cite{cl}. (We note, by the way, that \cite{kp} 
used moving planes to prove a symmetry result  on certain bounded domains in hyperbolic space.)  
\begin{theo}
Let $u$ be a global solution of (\ref{eq:heq}) on $\HH^n$.
\begin{itemize}
\item[i)] If $u$ is bounded and $\lim_{z \to \infty} u(z) = \eta$, where either $\eta = +1$ or $\eta = -1$,
then $u \equiv \eta$. 
\item[ii)] Suppose that $f$ is Lipschitz, with Lipschitz constant $L \leq (n-1)^2/4$. Suppose that $|u(z)|/(1-|z|)^{\alpha} \leq C$ 
for some $\alpha \in (\beta_-, \frac12(n-1))$. Then $u \equiv 0$. The same conclusion is valid if $f \in \calC^{1,\delta}$
for some $\delta > 0$ and $\alpha \in (\alpha_-,\frac12 (n-1))$.  Here $\alpha_-$ and $\beta_-$ are defined
in (\ref{eq:indroots}). 
\end{itemize} 
\label{th:ecs}
\end{theo}
\noindent{\bf Proof:} We prove ii) first. Consider $u$ in the ball model. Fix any hyperplane $P \subset 
\HH^n$ and let $R: \HH^n \to \HH^n$ be the reflection across $P$. Now define $u_R = u\circ R$, which
we consider as a function in one of the two components $\Omega$ of $\HH^n \setminus P$. 
Define $w_R = u_R - u$. We claim that $w_R \geq 0$ in $\Omega$. Interchanging the roles of the
two components, we see that if this claim holds, then $u$ is invariant with respect to reflections
about $P$, and since $P$ is arbitrary, we may conclude that $u$ is constant. 

Using that $f$ is Lipschitz, we calculate that
\[
\Delta w_R = f(u_R) - f(u) \leq L |w_R|.
\]
Suppose that $w_R < 0$ in some subdomain $\Omega' \subset \Omega$. Then $|w_R| = -w_R$ there,
hence $\Delta w_R + L w_R \leq 0$ in $\Omega'$.  Consider the function $v = w_R/\phi_\alpha$,
where $\phi_\alpha$ is the spherical function defined at the beginning of this section and 
$\alpha$ is as in the statement of this theorem. Then
\[
\nabla v = \frac{\nabla w_R}{\phi_\alpha} - \frac{w_R \nabla \phi_\alpha}{\phi_\alpha^2}, 
\]
\[
\Delta v  =  \frac{\Delta w_R}{\phi_\alpha} - 2 \frac{\nabla w_R \cdot \nabla \phi_\alpha}{\phi_\alpha^2}
- \frac{w_R \Delta \phi_\alpha}{\phi_\alpha^2} + 2 \frac{w_R |\nabla \phi_\alpha|^2}{\phi_\alpha^3},
\]
so, using the equations for $w_t$ and $\phi_\alpha$, 
\[
\Delta v + X \cdot \nabla v + (\alpha(\alpha - (n-1)) + L)v \leq 0,
\]
where $X = 2 \nabla \phi_\alpha/\phi_\alpha$. It is not hard to check, using the regularity of 
$\phi_\alpha$, that $X$ satisfies the hypothesis in Lemma~\ref{le:lemma2}. If $\alpha > \beta_-$, where $\beta_-$ 
is defined in (\ref{eq:indroots}), then $\alpha(\alpha - (n-1)) + L < 0$. Hence by this Lemma,
$v \geq 0$ in $\Omega'$, so $w_R$ is also nonnegative, which is what we claimed. 

To prove the final assertion of this theorem,  it suffices to assume that $\alpha > \alpha_-$ 
if $f$ is not just Lipschitz, but $\calC^{1,\delta}$ for some $\delta > 0$, we invoke a regularity theorem
from \cite{M-edge}, as implemented in \cite{M-rsyp}.  To set this up, write (\ref{eq:heq}) as a perturbation 
of the solution which is identically $0$:
\[
\Delta u - f'(0)u = Q(u).
\]
Let $\Lambda^{k,\nu}(\HH^n)$ denote the H\"older space of order $k+\nu$, $k \in {\mathbb N}$, $0 < \nu < 1$, where 
all derivatives and difference quotients are computed with respect to the hyperbolic metric. Denote by 
$\rho^s \Lambda^{k,\nu}$ the set of all functions $u = \rho^s v$ where $v \in \Lambda^{k,\nu}$; for convenience here and 
below, we denote by $\rho = 1-|z|^2$ a smooth defining function for $\del_\infty\HH^n$, and we shall also assume
that $\nu = \delta$ to reduce the number of different indices. The precise result we quote from these papers is as 
follows: fix $\lambda < (n-1)^2/4$, and let $\alpha_\pm$ be the two indicial roots for 
$\Delta + \lambda$ computed as in (\ref{eq:quadratics}) and (\ref{eq:indroots}); then for any $\alpha \in (\alpha_-,
\alpha_+)$, the mapping
\[
\Delta + \lambda: \rho^{\alpha} \Lambda^{2,\nu}(\HH^n) \longrightarrow \rho^{\alpha}\Lambda^{0,\nu}(\HH^n)
\]
is an isomorphism. To apply this, observe first that since $f \in \calC^{1,\nu}$, if $|u| \leq C(1-|z|)^s$ for 
any $s > 0$, then $|Q(u)| \leq C(1-|z|)^{(1+\nu)s}$. Now if $|u| \leq C\rho^s$ for $s > \alpha_-$ and
$\Delta u = f(u)$, then it follows from scale-invariant Schauder estimates (namely, just ordinary Schauder
estimates in balls of unit size with respect to hyperbolic distance) that $u \in \rho^s\Lambda^{3,\nu}(\HH^n)$
(i.e.\ precisely two derivatives more than the regularity of $f$). Now suppose $u \in \rho^\alpha \Lambda^{2,\nu}$ 
is a solution; then $Q(u) \in \rho^{\alpha(1+\nu)}\Lambda^{2,\nu}$, hence if $\alpha' := \alpha(1+\nu) < \alpha_+$,
we deduce that $u \in \rho^{\alpha'}\Lambda^{2,\nu}$. Continuing in this way a finite number of times implies
eventually that $u \in \rho^{\alpha''}\Lambda^{2,\nu}$ for some $\alpha'' > \beta_-$, and now we may use the previous argument. 

To prove i) we must argue in a slightly different way since we do not wish to use any information
about the size of the Lipschitz constant for $f$.  To be definite, fix $\eta = 1$; the argument
for $\eta = -1$ is identical. 

First fix any geodesic $\gamma$ and let $P_t$ denote the family of totally geodesic hyperplanes
perpendicular to $\gamma$. Let $R_t$ be the reflection across $P_t$, and set $u_t = u \circ R_t$.
Let $\Omega_t^+$ denote the `forward' component of $\HH^n \setminus P_t$, i.e.\ the component
swept out by the $P_s$, $s > t$, and $\Omega_t^-$ the other component. 
We claim that when $t \ll 0$, $u_t \geq u$ in $\Omega_t^+$.  To see this, note that on the one hand, 
for any given $\e > 0$, $1 > u \geq 1- \e$ when $t$ is sufficiently negative. Suppose now that
$u_t < u$ in some subdomain $\Omega' \subset \Omega_t^+$. Thus $u(z)$ also lies in $(1-\e,1)$ 
for $z \in \Omega'$.

The difference $w_t = u_t - u$ satisfies the equation 
\[
\Delta w_t = I(z) w_t, \qquad I(z) = \int_0^1 f'(u_t (z) + \tau (u(z) - u_t(z))\, d\tau.
\]
By assumption, $f(s)$ is strictly increasing in some neighbourhood $(1-\e,1)$, and the argument
of $f'$ lies in this interval when $z \in \Omega'$, so $I(z) \geq 0$ for $z \in \Omega'$, 
and hence by Lemma \ref{le:lemma2},  $w_t \geq 0$. 

Now let $t_0$ be the supremum of values in $\RR$ for which $u_t \geq u$ in $\Omega_t^+$. 
Clearly, unless $u$ is constant, $t_0$ must be finite since when $t \gg 0$, the roles
of the components $\Omega_t^\pm$ and $u_t$ and $u$ are interchanged and we would have the
reverse inequality.  Either $P_{t_0}$ must be a plane of symmetry for $u$ or else
$u_{t_0} > u$ in $\Omega_{t_0}^+$. To see that this latter case is impossible, we must
show that we can push $P_t$ a bit further while still preserving this inequality,
which contradicts that $t_0$ is the supremum. 
For this it suffices to show that the normal derivative of $w_{t_0}$ is bounded below
by a positive constant along all of $P_{t_0}$. If that were to fail, there would be
a sequence of points $q_j \in P_{t_0}$ such that $\del_{\nu} w_{t_0} \searrow 0$.
Choose a sequence of M\"obius transformations which preserve $P_{t_0}$ and which
move $q_j$ to a fixed point $q$. Applying Lemmas \ref{le:lemma1} and \ref{le:lemma2},
we obtain a limiting function $w$ which satisfies $\Delta w = I w$ with $w \geq 0$
and $\del_\nu w = 0$ at $q$, which violates the Hopf boundary point lemma. 

We have now shown that one of the hyperplanes $P_t$ orthogonal to the geodesic $\gamma$ must 
be a plane of symmetry for $u$. This is true for every geodesic $\gamma$, and we claim
that this implies that $u$ is in fact identically constant. 

The proof of this last fact is actually slightly easier than in Euclidean space since
in some sense there are `more hyperplanes' in hyperbolic space. As shown above, for every
(unoriented) geodesic $\gamma$ in $\HH^n$, the solution $u$ is symmetric around some 
hyperplane $P_\gamma$ orthogonal to $\gamma$. 
Suppose that two of these hyperplanes, $P_1$ and $P_2$, do not intersect at any finite
point of $\HH^n$ or on the asymptotic boundary $S^{n-1}$. Let $c$ be the unique geodesic 
perpendicular to both of them. The composition of the reflections across $P_1$ and
$P_2$ is a translation $h$ along $c$ which has the property that for any $z \in \HH^n$,
$h^k(z) \to \infty$. However, $u(h(z)) = u(z)$ for all $z$ and $u(h^k(z)) \to 1$, so
$u \equiv 1$. The same argument may be applied if $P_1$ and $P_2$ intersect at infinity,
but not at an interior point; then $h$ is a parabolic element, but it still has the
property that every orbit $\{h^k(z)\}$ is unbounded, so the conclusion is the same.
In fact, we see that the crucial property we need is that the group of motions $\Gamma$ of
$\HH^n$ generated by all the reflections $R_\gamma$ across the various planes of symmetry
$P_\gamma$ satisfies that the orbit of any point $z \in \HH^n$ under $\Gamma$ is unbounded.
Elementary hyperbolic geometry shows that this property is true for the group generated by
any three reflections $R_1$, $R_2$, $R_3$, provided the corresponding planes do not all pass 
through the same point.
In other words, either we can find some triplet of symmetry planes $P_1$, $P_2$, $P_3$
which do not pass through the same point, in which case we conclude that $u \equiv 1$,
or else all symmetry planes pass through some fixed point $0 \in \HH^n$. In this last
case, we conclude that $u$ is radial around this point. However, we showed in \S 2
that there are no nonconstant radial solutions to this equation, so once again, $u$
must be constant.  

This concludes the proof in all cases. 
\hfill $\Box$ 

\subsection{The parabolic case} 
We have already noted that the only situation in which it is possible to expect a
one-dimensional symmetry result when the boundary values are $H_p$-invariant is
when $u \to 0$ on $S^{n-1}\setminus \{q\}$ and $u \to \pm 1$ as $z$ tends nontangentially to $q$.  

\begin{theo} 
Suppose that $f$ is Lipschitz, with Lipschitz constant $L \leq (n-1)^2/4$. Let $\alpha > \alpha_-$,
and suppose that $u$ is a global 
bounded solution of (\ref{eq:heq}) satisfying $|u(z)|/\sigma_\alpha \leq C$ on all of $\HH^n$,
where the pole of $\HH^n$ is placed at some point $q \in S^{n-1}$, and moreover that $u(z) \to 1$
as $z \to q$ in any nontangential approach region in the ball. Then $u$ is $H_p$ invariant, 
hence in the upper half-space model where $q$ is at infinity, $u(x,y) = U(x)$ where $U$ is the 
function obtained in Proposition \ref{pr:pc}. 
\label{th:pcs}
\end{theo}
\noindent{\bf Proof:}  The proof is very close to that for part ii) of Theorem \ref{th:ecs}. 
Consider $u$ in the upper half-space model, with the exceptional point $q$ placed at infinity. 
Fix any vertical plane $P$ given by $\ell \cdot y = c$ and let $R$ denote the reflection across 
this plane and $u_R = u \circ R$. As before, we set $w_R=u_R-u$, so that by hypothesis 
$\Delta w_R= f(u_R) - f(u) \leq L |w_R|$. Let $\Omega^\pm$ be the two components of
$\HH^n \setminus P$. We claim that $w_R \geq 0$ in $\Omega^+$, say. 

The proof is very similar indeed to what we had done before. If $w_R < 0$ in some subdomain $\Omega'$,
then we introduce $v = w_R/\sigma_\alpha$; by the hypothesis, $|v| \leq C$, and $v=0$ on $\Omega' \cap P$,
and $v \geq -C$ in the entire subdomain. Since $\Delta v + X \cdot \nabla v \leq 0$ where $X$ is
a left-invariant vector field, we conclude that $v \geq 0$, a contradiction.   Applying the
same argument in the domain $\Omega^-$ shows that $u_R \equiv u$, so $P$ is a plane of symmetry.
Since every vertical plane is now a plane of symmetry, $u$ is independent of $y$, and hence
a function of $x$ alone.

It is very important here that $w_R \geq -C \sigma_\alpha$, so that we may control
the behaviour of $v$ as $(x,y) \to \infty$ nontangentially, e.g.\ as $y \to \infty$ with $x$ fixed. 
\hfill $\Box$. 

Note that even for ODE solutions, tangential limits of $u$ at $q$ can take on other values in $(-1,1)$. 

\subsection{Hyperbolic case}
The corresponding symmetry result in the hyperbolic case is less standard; the proof uses what could be 
called a `stretching method', in analogy to the `sliding method' in the Euclidean setting 

\begin{theo}
Let $u$ be a solution of (\ref{eq:heq}) satisfying the asymptotic boundary condition 
that $u(z)$ converges to $\pm 1$ at any point in the interior of the two
hemispheres $S^\pm \subset \del \HH^n$. Let $\HH^{n-1}$ be the totally geodesic
subspace with boundary $S_- \cap S_+$, and $t$ the signed distance function to
this subspace in $\HH^n$. Suppose that $f$ is a nonnegative $\calC^1$ 
function on $\RR$ such that $f(1)=f(-1)=0$ and $f'(1)>0$, $f'(-1)>0$. 
Then $u(z)$ depends only on $t$, i.e.\ $u(z) = U(t(z))$, where $U$ 
satisfies the ordinary differential equation 
\[
U''(t) + (n-1)\frac{\sinh t}{\cosh t} U'(t) = f(U), \qquad \lim_{t \to \pm \infty} U(t) = \pm 1.
\]
\end{theo}
 
\noindent {\bf Proof:} Let $P_{-} \in S_-$ and $P_+ \in S_+$ be arbitrary points,
and suppose that $\gamma(t)$ is any curve of constant geodesic curvature in $\HH^n$ 
such that $\lim_{t \to \pm \infty} \gamma(t) = P_{\pm}$. We shall first prove that 
$u\circ \gamma(t)$ is (weakly) monotone increasing in $t$. The conclusion is not
difficult to prove from this.

For this step, we work in the upper half-space model. Compose with a M\"obius transformation so 
that $P_- = 0$ and $P_+ = \infty$; the images of
$S_-$ and $S_+$ (which we denote by the same symbols) are then some ball in $\RR^{n-1}$ containing 
$0$, but not necessarily at its center, and the exterior of this ball (union $\infty$), respectively, 
and the curve $\gamma$ is transformed into some ray in the upper half-space emanating from $0$.

Fix $\delta > 0$ such that $f'(s)>0$ when $s \in [-1,-1+\delta)\cup (1-\delta,1]$, and choose $A \in (0,1)$
so that $u(z)<-1+\delta$ for $z \in D_-(A) := \{z \in \RR^n_+: |z| < A\}$ and $u(z) >1-\delta$ for 
$z \in D_+(A) := \{z \in \RR^n_+: |z| > 1/A\}$. We let $\del_i D_\pm (A)$ denote that portion of 
$\del D_\pm (A)$ lying in $\HH^n$. 

Now choose any $R > 1/A^2$ and define $u_R(z) = u(Rz)$. We claim that $u_R(z)\geq u(z)$ 
for all $z \in \HH^n$.  We verify this in $D_-(A)$ and $\HH^n \setminus D_-(A)$ separately.
First note that if $z \in \del D_-(A)$, then $u_R(z) > 1-\delta > -1 + \delta > u(z)$. 
Next, suppose that the subset $\Omega = \{z \in D_-(A): u_R(z) < u(z)\} \subset D_-(A)$ 
is nonempty, and define $w_R = u_R - u$. Then in $\Omega$, 
\begin{equation}
\Delta w_R = I(z)w_R, \qquad I(z) = \int_0^1f'(u(z)+t(u_R(z)-u(z)))\,dt,
\label{eq:uru}
\end{equation}
while $w_R=0$ on the entire boundary of $\Omega$, including the portion lying in $\del \HH^n$. 
Furthermore, if $z \in \Omega$, then $u_R(z) < u(z) < -1 + \delta$, and hence $u(z) + t(u_R(z) - u(z)) < 
-1 + \delta$ for all $0 \leq t \leq 1$. Since $\delta < 1$, $I(z)\geq 0$ on $\Omega$.  
We can now apply the maximum principle to (\ref{eq:uru}) to conclude that $u_R - u \geq 0$ in 
$\Omega$, which is a contradiction to the definition of this set, and hence $\Omega=\emptyset$.  
Finally, if $z \in \HH^n \setminus D_-(A)$, then $u_R(z) >1-\delta$, and a similar argument proves 
that $u_R\geq u$ holds outside of $D_-(A)$ as well. This proves the claim. 

Now define $\bar R=\inf\{R>1: u_R(z) \geq u(z) \ \forall\, z \in \HH^n\}$. We have just
established that this set is nonempty, so that $\bar R \in [1,\infty)$ is well defined. 
We claim that $\bar R=1$. 

Suppose by way of contradiction that $\bar R > 1$.  Define $\bar u = u_{\bar R}$, $w_R = u_R - u$
and $\bar w = \bar u - u$. By continuity, $\bar w \geq 0$ on $\HH^n$, and $\bar w$ satisfies 
an elliptic equation, as in (\ref{eq:uru}), with coefficient $I(z)$. Decomposing $I = I_+ - I_-$,
we see that $\Delta \bar w = (I_+ - I_-)\bar w \leq I_+ \bar w$, so by the strict minimum
principle for supersolutions, $\bar w >0$ on $\HH^n$.

If $\bar R>1$, there would exist a sequence $R_k \nearrow \bar R$ and points $z_k \in \HH^n$ so that
the function $w_{R_k} := w_k$ attains a strictly negative absolute minimum value at $z_k$.  If some
subsequence of the points $z_k$ were to converge to a point $\bar z \in \HH^n$, then we would
have $\bar w (\bar z) =0$, which we have shown cannot happen. Thus $z_k$ diverges in $\HH^n$.

To renormalize, choose a sequence of M\"obius transformations, $T_k$, so that $T_k(z_k)$ is 
equal to some fixed point $\bar z \in \HH^n$. Denoting $w_k \circ T_k$ by $\tw_k$, then as before, 
each $\tw_k$ satisfies an equation $\Delta \tw_k = \tI_k(z) \tw_k$, where $\tI_k = I_k \circ T_k$;
furthermore, $\inf \tw_k = \tw_k(\bar z) = -\epsilon_k \nearrow 0$. 

Using Lemmas \ref{le:lemma1} and \ref{le:lemma2} again, we may assume that 
$\tu_k \to \tu$, $\tw_k \to \tw$ and $\tI_k \to \tI$, and
\[
\Delta \tw= \tI \tw, \qquad \tw \geq 0, \qquad \tw(\bar z) = 0.
\]
Arguing as above, we may apply the strict minimum principle for supersolutions again to conclude that $\tw \equiv 0$. 
This shows that $\tI = f'(\tu)$; however, since this function
does not necesssarily have a sign, we must proceed further. In fact, we derive a contradiction
by considering the function $W_k = \epsilon_k^{-1} \tw_k$ (where $-\epsilon_k = \inf \tw_k$). Clearly 
$\Delta W_k =I_k W_k$, but now $\min\, W_k= -1$. As before, we can extract a subsequence converging
in $\calC^{2,\alpha}$ on every compact set to some function $\bar{W}$, where
\begin{equation}
\Delta \bar W= \tI \bar W, \qquad \inf \bar{W} = \bar W(\bar z) = -1. 
\label{eq:maxbW}
\end{equation}

Now let us show why such a limit cannot exist.  There are two cases to consider, depending on
the limiting behaviour of the sequence $z_k$ in $\overline{\HH^n}\setminus \HH^n$:
\begin{itemize}
\item[i)] Suppose that $z_k$ converges to a point $z_\infty$ which lies in either of the open half-balls
$S_\pm$ where $u$ is continuous and equals $\pm 1$. In this case, $\tI(z) = f'( \lim \tu_k(z)) = f'(\pm 1) > 0$,  
so (\ref{eq:maxbW}) is impossible. 
\item[ii)] Suppose on the other hand that $z_k$ converges to some point $z_\infty$ in the intersection of
the two hemispheres $S_+ \cap S_-$. Now recall that (using upper half-space notation again), since $R_k \geq c > 1$
for all $k$, then $R_k z_k$ must converge to a point in $S_+$, so $u(R_k z_k) \to 1$. However, by definition,
$u(R_kz_k) - u(z_k) < 0$, so $u(z_k) \to 1$ as well. (This is not incompatible with $z_\infty$ lying
on $S_+ \cap S_-$, but rather, indicates that it approaches this intersection tangent to the $S_+$ side.)
Thus here too we conclude that $\tI(z) > 0$, at least in a neighbourhood of $\bar z$, so
(\ref{eq:maxbW}) is again impossible. 
\end{itemize}

We have now ruled out the possibility that $\bar R > 1$, so $\bar R = 1$. In other words, we have proved that 
$u_R\geq u$ for all $R>1$; hence
\[
\left. \frac{d\,}{dR} \right|_{R=1} u(Rz) = r\del_r u(z) \geq 0
\]
for every $z \in \HH^n$; here $r$ is the Euclidean radial upper-half-space variable. This shows that
$u$ is nondecreasing along any ray emanating from the origin, i.e.\ along any curve of constant
geodesic curvature joining $P_-$ to $P_+$. 

We now prove that $u$ is a function of only one variable. Working again in the ball model,
suppose that $S_-$ and $S_+$ are the lower and upper hemispheres of the boundary, respectively. Let $P$
be any point in the interior of the ball, and let $\Sigma$ be the spherical cap which passes through $P$ 
and $S_- \cap S_+$. Thus $\Sigma$ is a hypersurface equidistant from the totally geodesic
copy of $\HH^{n-1}$ which has the same asymptotic boundary $S_+ \cap S_-$. Let $\Pi$ be any two-dimensional
plane passing through the origin of the ball and the point $P$. Then $\Pi \cap \Sigma$ is a curve
of constant geodesic curvature $\gamma$ in $\HH^n$ passing through $P$ and limiting on two points
$Q, Q' \in S_- \cap S_+$. It is easy to see geometrically (for example by tilting $\Sigma$ slightly,
but still ensuring that it passes through $P$), that we can approximate $\gamma$ by two sequences
of curves of constant geodesic curvature, $\gamma^-_j(t)$ and $\gamma^+_j(t)$ such that
$\gamma^\pm_j(0) = P$ for all $j$, 
\[
\lim_{t \to -\infty} \gamma_j^-(t), \lim_{t \to +\infty} \gamma_j^+(t) \in S_-, \qquad
\lim_{t \to +\infty} \gamma_j^-(t), \lim_{t \to -\infty} \gamma_j^+(t) \in S_+,
\]
and finally 
\[
T_P \Sigma \ni X = \lim_{j \to \infty} (\gamma_j^-)'(0) = - \lim_{j \to \infty} (\gamma_j^+)'(0).
\]
Since $u \circ \gamma_j^-(t)$ and $u \circ \gamma_j^+(t)$ are both nondecreasing, we see that
$\nabla u_P \cdot X = 0$. However, $X$ can be chosen arbitrarily in $T_P \Sigma$, which shows
that $\nabla u(P)$ is orthogonal to $\Sigma$.

We have now proved that if $\Sigma_t$ is the foliation of $\HH^n$ by hypersurfaces which
are of (signed) distance $t$ from the totally geodesic copy of $\HH^{n-1}$ with boundary
$S_- \cap S_+$, then each $\Sigma_t$ is a level set of $u$. In other words, $u$ is a function
of the distance $t$ alone. This concludes the entire proof. \hfill $\Box$ 

\section{Non-symmetric solutions}
We have already indicated that the decay hypotheses in parts i) and ii) of Theorem 1.2 are necessary
since in fact there is a very large class of non-symmetric solutions of (\ref{eq:heq}) the elements of
which satisfy the bounds 
\[
\begin{array}{rcl}
&a)\ &|u(z)| \leq C (1-|z|^2)^{\alpha_-},  \\
&b)\ &|u(z)| \leq C x^{\alpha_-}, \qquad \mbox{and}\qquad \lim_{x \to \infty} u(x,y) = 1
\end{array}
\]
respectively. In this section we give a brief sketch of the construction of these solutions
under the slightly stronger assumption that $f \in \calC^{1,\delta}$ (this is probably not necessary 
but is used in the proof). 

This construction is a simple perturbation argument using the implicit function theorem.
Suppose that $u_0$ is a known solution of (\ref{eq:heq}). Then any nearby solution $u$
may be decomposed as $u = u_0 + v$ where $v$ is small, both globally and perhaps also in
terms of its asymptotic decay rate. Expanding $f$ into its first order Taylor polynomial
and an error term, the equation $\Delta u = f(u)$ can be rewritten as
\[
\Delta v = f'(u_0)v + Q(u_0;v),
\]
where $Q$ at least has the property that if $|v| \leq \e$ then $|Q(v)| \leq C\e^{1+\delta}$,
and the notation $v \mapsto Q(u_0;v)$ is meant to indicate the dependence of $Q$ on $u_0$. 

We shall actually use a slightly finer decomposition. Let $\phi$ any element of the nullspace
of $\Delta - f'(u_0)$; in the examples below, $\phi$ will have the same asymptotic behaviour
as $u_0$, but will not decay faster. Now write $u = u_0 + \phi + w$, where $w$ lies in some 
space of functions which decay faster than $u_0 + \phi$ at infinity. We rewrite the equation
one last time as
\begin{equation}
(\Delta - f'(u_0))w = Q(u_0;\phi + w).
\label{eq:perturb}
\end{equation}

The existence of nearby solutions can therefore be deduced from the mapping properties of
the operator $L = \Delta - f'(u_0)$ and from the asymptotic properties of solutions of 
$L\phi = 0$. We now describe these in the two cases of interest. The reason that this
sort of analysis works is because we consider perturbations of solutions $u_0$ which
are asymptotically equal to $0$ on most or all of $\del_\infty \HH^n$ and because $0$
is an unstable critical point of the potential function $F$. We shall use the weighted
H\"older spaces $\rho^s \Lambda^{k,\nu}$ introduced in \S 3.1. 

For the elliptic case, let $u_0 \equiv 0$. Recall the Poisson transform $P$ on
hyperbolic space. It is well-known that $P$ provides a correspondence 
between $\calC^{2,\nu}(S^{n-1})$ and the space of solutions of $(\Delta - f'(0))\phi = 0$
with $\phi \sim \phi_0 \rho^{\alpha_-}$, $\phi_0 \in \calC^{2,\mu}(S^{n-1})$. 

\begin{prop}
Let $u_0 \equiv 0$ and fix $\alpha \in (\alpha_-,\alpha_+)$, where $\alpha_\pm$ are
described in (\ref{eq:quadratics}) and (\ref{eq:indroots}). There is a neighbourhood 
of $0$ $\calU \subset \calC^{2,\mu}(S^{n-1})$ and a map $G: \calU \to 
\calV \subset \rho^\alpha \Lambda^{2,\nu}(\HH^n)$ so that all solutions $u \in
\rho^{\alpha_-}\Lambda^{2,\nu}$ to (\ref{eq:perturb}) close to $0$ are of the form 
$u = P(\phi_0) + G(\phi_0)$. 
\end{prop}
\noindent{\bf Proof:} First note that $(\phi_0,w) \mapsto Q(0;P(\phi_0) + w)$
is a continuous mapping from $\calU \times \calV$ to $\rho^{\alpha}\Lambda^{2,\nu}(\HH^n)$
provided $\alpha \leq (1+\delta)\alpha_-$. (Here $\calU$ and $\calV$ are small balls
in their respective function spaces, as in the statement of this result.) 
Furthermore, as described in \S 3.1 and proved in \cite{M-edge}, since $\lambda = 
-f'(0) < (n-1)^2/4$, 
\[
\Delta + \lambda: \rho^\alpha \Lambda^{2,\nu}(\HH^n) \longrightarrow \rho^\alpha \Lambda^{0,\nu}(\HH^n)
\]
is an isomorphism. We now rewrite (\ref{eq:perturb}) as
\[
w = (\Delta + \lambda)^{-1}Q(0;P(\phi_0) + w)
\]
and, for each $\phi_0$ sufficiently close to $0$, obtain a solution by showing that the map
on the right is a contraction mapping. \hfill $\Box$. 

\medskip

In the parabolic case we follow a similar strategy. Since we are not trying to obtain
the most general possible result, but only wish to exhibit the abundance of possible solutions,
we simplify the analysis by restricting the allowable perturbations to have a very special form.
Let $u_0(x)$ denote the ODE solution obtained in Theorem 1.1, part ii), and consider
nearby solutions $u(x,y)$ which are invariant with respect to the integer lattice
$\Gamma = {\mathbb Z}^{n-1}$ acting on $y \in \RR^{n-1}$. In other words, we consider 
the problem on the quotient $M = \HH^n / \Gamma$. This manifold $M$ is diffeomorphic
to the product $\RR^+ \times \TT$, $\TT = T^{n-1} = \RR^{n-1}/{\mathbb Z}^{n-1}$; it has two ends,
one a `large' conformally compact end to which the same sort of analysis as above applies, 
and the other a finite volume hyperbolic cusp end. As a further simplification, we shall
consider perturbations of the form $u(x,y) = u_0(x) + v(x,y)$ where $v$ is not only periodic
in $y$, hence descends to this quotient, but in addition, $\int_{\TT} v(x,y)\, dy = 0$ for all $x > 0$.
(This is sometimes called the cusp form condition.) 

Let $\phi_0 \in \calC^{2,\alpha}(\TT)$ be any function with vanishing integral. There is once
again a Poisson transform $P$ which maps $\phi_0$ to the solution $\phi(x,y)$ to
$(\Delta - f'(u_0))\phi = 0$ with $\phi(x,y) \sim \phi_0(y)x^{\alpha_-}$ as $x \searrow 0$.
This solution $\phi$ satisfies the cusp form condition, and decays exponentially as $x \nearrow \infty$.
All of these facts may be proved by separating variables and using standard ODE arguments.

The next step is to construct a parametrix $G$ for $L := \Delta - f'(u_0)$ by pasting together
parametrices for $\Delta - f'(0)$ near $x=0$ and for $\Delta - f'(1)$
for $x \gg 0$. Using the same sort of analysis from \cite{M-edge} as employed in
the previous result, there is an inverse $G_0$ for $\Delta - f'(0)$ on all of $M$;
similarly, there is also an inverse $G_\infty$ for $\Delta - f'(1)$ since $f'(1)>0$.
Now choose a partition of unity $\chi_0 + \chi_\infty = 1$ where $\chi_0 = 0$ for
$x \gg 0$ and $\chi_\infty = 0 $ for $x \leq 1$, say. Also choose functions
$\tilde{\chi}_0$, and $\tilde{\chi}_\infty$ such that $\tilde{\chi}_j = 1$ on the
support of $\chi_j$, $j=0,\infty$, and $\tilde{\chi}_0$ vanishes for $x \gg 0$, $\tilde{\chi}_\infty = 0$
for $x \leq 1/2$, say. Then
\[
G = \tilde{\chi}_0 G_0 \chi_0 + \tilde{\chi}_\infty G_\infty \chi_\infty
\]
is a parametrix for $L$ in the sense that $L G = I - R$ where $R$ is
a smoothing operator which induces extra decay. We can choose all of these operators
to preserve the subspace of functions on $M$ which satisfy the cusp-form condition. 

To be more precise about mapping properties, set $\langle x \rangle = \sqrt{1+x^2}$, and
consider weighted spaces of the form $(x/\langle x \rangle)^\alpha \langle x \rangle^\beta \calC^{k,\nu}(M)$.
These functions grow or decay like $x^\alpha$ and $x^\beta$ as $x \searrow 0$, $x \nearrow \infty$,
respectively. Then
\[
G: (x/\langle x \rangle)^\alpha \langle x \rangle^\beta \Lambda^{0,\alpha}
\longrightarrow 
(x/\langle x \rangle)^\alpha \langle x \rangle^\beta \Lambda^{2,\alpha}
\]
is bounded for $\alpha \in (\alpha_-,\alpha_+)$ and $\beta \in (-\sqrt{f'(1)}, \sqrt{f'(1)})$.
Restricting to the subset of functions satisfying the cusp-form condition, we could even let
$\beta$ be an arbitrary real number. Furthermore, 
\[
R: (x/\langle x \rangle)^\alpha \langle x \rangle^\beta \Lambda^{0,\alpha}_0
\longrightarrow 
(x/\langle x \rangle)^{\alpha+1} \langle x \rangle^{\beta-1} \Lambda^{3,\alpha}_0,
\]
from which its compactness follows easily.  All of this shows that 
\[
L: (x/\langle x \rangle)^\alpha \langle x \rangle^\beta \Lambda^{2,\alpha}_0
\longrightarrow 
(x/\langle x \rangle)^{\alpha} \langle x \rangle^{\beta} \Lambda^{0,\alpha}_0,
\]
is Fredholm. It remains to show that this map is actually an isomorphism, 
which is true and can again be verified using separation of variables and ODE
arguments. 

The rest of the argument is exactly the same as in the previous case. This proves the
\begin{prop}
There is an infinite dimensional space of solutions of $\Delta u = f(u)$ on
$\HH^n / \Gamma$ which are small perturbations of the ODE solution $u_0(x)$. 
Each such solution has an asymptotic expansion of the form $u(x,y) \sim (c + 
\phi_0(y))x^{\alpha_-}$ as $x \searrow 0$, for all sufficiently small $\phi_0$
satisfying the cusp-form condition. 
\end{prop}

\end{document}